\documentclass[11pt,a4paper]{article}
\usepackage[utf8]{inputenc}
\usepackage[T1]{fontenc}
\usepackage{amsmath}
\usepackage{amssymb}
\usepackage{amsthm}

\usepackage{color}
\definecolor{DarkOlive}{rgb}{0.1047,0.2412,0.0064}
\definecolor{FireBrick}{rgb}{0.5812,0.0074,0.0083}
\definecolor{RoyalBlue}{rgb}{0.0236,0.0894,0.6179}
\definecolor{RoyalGreen}{rgb}{0.0236,0.6179,0.0894}
\definecolor{RoyalRed}{rgb}{0.6179,0.0236,0.0894}
\definecolor{LightBlue}{rgb}{0.8544,0.9511,1.0000}
\definecolor{Black}{rgb}{0.0,0.0,0.0}
\definecolor{MidnightBlue}{rgb}{0.0035,0.0020,0.1363}
\definecolor{FireBrick3}{rgb}{0.5812,0.0074,0.0083}
\definecolor{FireBrick4}{rgb}{0.2156,0.0023,0.0035}
\definecolor{Blue2}{rgb}{0.0000,0.0000,0.8644}
\definecolor{Navy}{rgb}{0.0000,0.0000,0.1927}
\definecolor{MediumBlue}{rgb}{0.0000,0.0000,0.6179}
\usepackage[
        bookmarks=false,pdftitle={Restrictions to a principal A1-subgroup},
        colorlinks=true,backref=false,breaklinks=true,linkcolor=MediumBlue,
        citecolor=FireBrick3,filecolor=RoyalRed,
        urlcolor=Blue2]{hyperref}

\theoremstyle{plain}
\newtheorem{Thm}{Theorem}
\newtheorem{Cor}[Thm]{Corollary}
\newtheorem{Pro}[Thm]{Proposition}
\newtheorem{La}[Thm]{Lemma}

\theoremstyle{definition}
\newtheorem{Def}[Thm]{Definition}
\newtheorem{Rem}[Thm]{Remark}

\newenvironment{Prf}{\noindent\textbf{Proof.}}%
  {\unskip\nobreak\hfil\penalty50\hskip2em\hbox{}\nobreak\hfil$\Box$%
   {\parfillskip=0pt \finalhyphendemerits=0 \par}\medskip}

\newcommand{\mset}[1]{\{\!\{ #1 \}\!\}}
\newcommand{\hgt}{ \textrm{ht} }

\newcommand{\Z}{\mathbb{Z}}
\newcommand{\R}{\mathbb{R}}
\newcommand{\C}{\mathbb{C}}
\newcommand{\Q}{\mathbb{Q}}

\newcommand{\SL}{\textrm{SL}}
\newcommand{\Gm}{\mathbb{G}_m}

\begin{document}
\title{Restrictions to a principal $A_1$-subgroup}
\author{Frank Lübeck\thanks{This is a contribution to 
Project-ID 286237555 -- TRR 195 -- by the
German Research Foundation}}
\maketitle
\begin{abstract}
We consider simple simply-connected algebraic groups $G$ over 
the complex numbers, and a principal $A_1$-subgroup $A \leq G$.
We show for which types of $G$ all its irreducible representations
$V(\lambda)$ ($\lambda$ a dominant weight) are determined by 
their restrictions to $A$ up to graph automorphisms. See Theorem~\ref{thmMain} below for 
the precise statement.
\end{abstract}

\section{Introduction}\label{secIntro}

Let $G$ be a connected reductive algebraic group defined over the complex
numbers, and assume that $G$ is simple and of simply-connected type. 
The group $G$ contains a \emph{principal $A_1$-subgroup}, that is, 
a closed subgroup which is the image of a homomorphism $\SL_2 \to G$ and
which contains a regular unipotent element of $G$. Such subgroups exist and
form a single conjugacy class of subgroups of 
$G$~\cite{jacobson,morozov,kostant}.

J.-P.~Serre~\cite{serreA1}, referring to~\cite{npp}, asked  
how much information about a
dominant weight $\lambda$ of $G$ is determined by the restriction of the
irreducible module $V(\lambda)$ of $G$ to a principal $A_1$-subgroup. 
In the article~\cite{npp} the authors stated
without proof that for groups of rank at most $4$ the restriction 
determines $\lambda$ up to graph automorphisms.
Here we prove the following theorem.

\begin{Thm}\label{thmMain}
For $G$ of rank at most $4$, or of type $B_l$ or $C_l$ for any $l$,
or of type $D_l$ with $l$ not a power of $4$,
or of exceptional type $E_6$, $E_7$ or $E_8$, 
let $\lambda$, $\mu$ be dominant weights of $G$. If the restrictions of
$V(\lambda)$ and $V(\mu)$ to a principal $A_1$-subgroup of $G$ are 
isomorphic then $\lambda = \mu$, or $G$ is of type $A_l$ ($l \leq 4$), 
$D_l$ or $E_6$ and $\mu$ is the image of $\lambda$ under a non-trivial graph 
automorphism.

The same statement is false for $G$ of type $A_l$ with $l \geq 5$,
and for $G$ of type $D_l$ with $l = 4^m$, $m \geq 2$.
\end{Thm}

By the results in Section~\ref{secRestr} this theorem is equivalent to a
purely combinatorial statement about the coroots of $G$, see 
Theorem~\ref{thmMultiset}. We prove it for exceptional groups in 
Section~\ref{secResults} (only a partial result for type $E_7$).
In Section~\ref{secClassical} we show with the same elementary methods
that the restriction also determines $\lambda$ up to graph automorphisms
for all simple groups of classical type and rank at most $4$, and we
describe counterexamples in type $A_l$, $l \geq 5$. In
Section~\ref{secInvariants} we give a proof for $G$ of type $E_7$, of type
$B_l$, $C_l$ for all $l$, and of type $D_l$ for $l$ not a power of $4$, using
invariant theory, and we construct counterexamples for type $D_l$ with
$l = 4^m$, $m \geq 2$.

\subsection{Notation}\label{secNotation}

We recall some facts and notation for which we refer to books on linear
algebraic groups, for example~\cite{spr,hum,milne,maltes}.

Let $B \subseteq G$ be a Borel subgroup and $T \subseteq B$ be a maximal
torus of $G$. Let $X = \textrm{Hom}(T, \Gm)$ be the character group and
$Y = \textrm{Hom}(\Gm, T)$ be the cocharacter group of $T$, both written
additively. For $x \in X$ and $y \in Y$ the composition $x \circ y$ is a
homomorphism $\Gm \to \Gm$, so there is a unique integer 
$\langle x, y \rangle \in \Z$ with
\[ x(y(t)) = t^{\langle x,y \rangle} \quad \textrm{for all } t \in \Gm. \]
This defines a non-degenerate pairing $\langle \cdot, \cdot \rangle
: X \times Y \to \Z$. To $T \subseteq G$ there is associated a root datum
$(X,\Phi,Y,\Phi^\vee)$, where $\Phi \subseteq X$ is the set of roots and
$\Phi^\vee \subseteq Y$ is the set of coroots of $G$, and there is a
bijection $\Phi \to \Phi^\vee$, $\alpha \mapsto \alpha^\vee$. 
The Weyl group $W = N_G(T)/T$ acts on $X$ and $Y$, the pairing is
$W$-invariant, and we write $\det(w) = \pm 1$ for the determinant of 
$w \in W$ acting on $X$.

Let $l$ be the rank of $G$. The Borel subgroup $B$ determines a set of
simple roots $\{\alpha_1, \ldots, \alpha_l\} \subseteq \Phi$; we
use the labelling of~\cite[Planches]{bou}. The corresponding coroots
$\{\alpha_1^\vee, \ldots, \alpha_l^\vee\}$ are a set of simple coroots.
Let $\Phi^+$ be the set of positive roots. Then 
$\Phi^\vee_+ = \{\alpha^\vee \mid \alpha \in \Phi^+\}$ is the set of 
coroots of the form $\alpha^\vee = b_1 \alpha_1^\vee + \ldots + b_l
\alpha_l^\vee$ with all $b_i \in \Z_{\geq 0}$, and $\Phi^\vee = \Phi^\vee_+
\cup -\Phi^\vee_+$.

Our assumption that $G$ is simply-connected means that the simple coroots
are a $\Z$-basis of $Y$. The elements $\{\omega_1,\ldots, \omega_l\}$ of
the dual basis of $X$ (with respect to $\langle \cdot, \cdot \rangle$) are
the fundamental weights. For a weight $\lambda = a_1 \omega_1 + \ldots +
a_l \omega_l \in X$ and $\alpha^\vee$ as above we have
\[ \langle \lambda, \alpha^\vee\rangle = a_1 b_1 + \ldots + a_l b_l. \]
The irreducible rational representations of $G$ are parameterized by the
set of dominant weights $\{\lambda = a_1 \omega_1 + \ldots + a_l \omega_l
\mid\; a_i \in \Z_{\geq 0}\}$. We write $V(\lambda)$ for the irreducible 
module with highest weight $\lambda$.

Let $\rho = \omega_1 + \ldots + \omega_l$; this is half the sum of the
positive roots. Similarly, let 
\[ 2\rho^\vee = \sum_{\alpha \in \Phi^+} \alpha^\vee \in Y. \] 
Applying~\cite[Ch.~VI, \S1.10]{bou} to the dual root system $\Phi^\vee$ 
we get $\langle \alpha_i, 2\rho^\vee \rangle = 2$ for $1 \leq i \leq l$.
Note that $\rho^\vee$ itself is in general not an element of $Y$.
For a positive coroot $\alpha^\vee = \sum_i b_i \alpha_i^\vee$ the integer
$\langle \rho, \alpha^\vee \rangle = b_1 + \ldots + b_l$ is the height of
$\alpha^\vee$ in $\Phi^\vee$.

\subsection{The principal $A_1$-subgroup}\label{secPrincipal}

There is a homomorphism $\phi: \SL_2 \to G$ such that 
\[ \phi\left(\begin{pmatrix} t & 0 \\ 0 & t^{-1}\end{pmatrix}\right) 
   = 2\rho^\vee(t) \quad \textrm{for all } t \in \Gm, \]
and its image is a principal $A_1$-subgroup of 
$G$~\cite{siebenthal,dynkin,kostant}, see also~\cite[\S2.2]{serre}.
Since all principal $A_1$-subgroups are conjugate in $G$, it follows that
every principal $A_1$-subgroup $A$ of $G$ has a maximal torus which is
conjugate to the image of $2\rho^\vee: \Gm \to T$.

The finite dimensional rational representations of $\SL_2$ are semisimple,
and they are determined up to isomorphism by their characters on the
diagonal torus. Since $A$ is an image of $\SL_2$, the same holds for
representations of~$A$. Therefore, the restriction of $V(\lambda)$ to $A$ is
determined up to isomorphism by the Laurent polynomial
\[ \chi_\lambda(2\rho^\vee(t)) \in \Z[t, t^{-1}], \]
where $\chi_\lambda$ denotes the character of $V(\lambda)$, and vice
versa. 

\begin{Rem}\label{remCharp}
Let $k$ be an algebraically closed field of characteristic $p>0$, and let
$G$ now be defined over $k$. If $p$ is at least the Coxeter number $h$ of $G$
then principal $A_1$-subgroups (closed subgroups of type $A_1$ containing a
regular unipotent element) exist and form a single conjugacy class of
subgroups of $G$~\cite{testerman,seitz}, see also~\cite{rt}. As above,
their maximal tori are conjugate to the image of $2\rho^\vee$. The
composition factors of the restriction of a $G$-module to a principal
$A_1$-subgroup are determined by the restriction of its character to the
image of $2\rho^\vee$. The Weyl modules of $G$ have the same characters as
the modules $V(\lambda)$ in characteristic $0$, 
see~\cite[II.5.10]{jantzen}. Therefore, all results in this note
also hold for the composition factors of the restrictions of Weyl modules
to principal $A_1$-subgroups in characteristic $p \geq h$.
\end{Rem}

\section{Description of restrictions}\label{secRestr}

For a rational $T$-module $M$ with weight spaces $M_x$, $x \in X$, we
consider its character as an element of the group ring $\Z[X]$ with basis 
$\{e^x \mid x \in X\}$: $\textrm{ch}(M) = \sum_{x \in X} \dim(M_x)\, e^x$.
Evaluating characters at the elements $2\rho^\vee(t)$, $t \in \Gm$, 
corresponds to applying the ring homomorphism
\[ \psi: \Z[X] \to \Z[t,t^{-1}],\quad e^x \mapsto t^{\langle x, 2\rho^\vee
\rangle}. \]
For $\theta \in X$ we set $A_\theta = \sum_{w \in W} \det(w)\, e^{w\theta}
\in \Z[X]$. Weyl's character formula, see~\cite[24.3]{humlie} 
or~\cite[II.5.10]{jantzen}, says that for a dominant weight $\lambda$ and
its irreducible module $V(\lambda)$ with character $\chi_\lambda$ we have
\[ A_\rho \cdot \chi_\lambda = A_{\lambda + \rho}. \]

\begin{La}[{\cite[2.6]{serre}}]\label{laDenom}
For $\theta \in X$ we have
\[ \psi(A_\theta) = \prod_{\alpha \in \Phi^+} \left(t^{\langle \theta,
\alpha^\vee \rangle} - t^{-\langle \theta, \alpha^\vee \rangle}\right). \]
\end{La}
\begin{Prf}
The coroots $\Phi^\vee$ form a root system in $\R \otimes Y$ with Weyl group
$W$, set of positive roots $\Phi^\vee_+$ and half sum of positive roots
$\rho^\vee$. The Weyl denominator identity, see~\cite[24.1(A) and
24.3]{humlie}, applied to this root system yields
\[ \sum_{w \in W} \det(w)\, e^{w\rho^\vee} = \prod_{\alpha \in \Phi^+}
\left(e^{\alpha^\vee/2} - e^{-\alpha^\vee/2}\right) \]
in the group ring $\Z[\frac{1}{2}Y]$. Now we apply the ring homomorphism
$\Z[\frac{1}{2}Y] \to \Z[t,t^{-1}]$, $e^y \mapsto t^{2\langle \theta, y
\rangle}$. The right hand side maps to the claimed product, and the left
hand side maps to 
\[ \sum_{w \in W} \det(w)\, t^{\langle \theta, w(2\rho^\vee)\rangle} = 
\sum_{w \in W} \det(w^{-1})\, t^{\langle w^{-1}\theta, 2\rho^\vee\rangle} =
\psi(A_\theta). \]
\end{Prf}

\begin{Thm}[{\cite[Thm.~2.8]{serre}}]\label{thmSerre}
Let $\lambda$ be a dominant weight. Then
\[ \chi_\lambda(2\rho^\vee(t)) = 
\prod_{\alpha \in \Phi^+} 
\frac{t^{\langle \lambda+\rho, \alpha^\vee \rangle} - 
      t^{-\langle \lambda+\rho, \alpha^\vee \rangle}}
     {t^{\langle \rho, \alpha^\vee \rangle} - 
      t^{-\langle \rho, \alpha^\vee \rangle}}
= t^{-\langle \lambda, 2\rho^\vee \rangle}
\prod_{\alpha \in \Phi^+} 
\frac{t^{2\langle \lambda+\rho, \alpha^\vee \rangle} - 1}
     {t^{2\langle \rho, \alpha^\vee \rangle} - 1}. \]
\end{Thm}
\begin{Prf}
We apply $\psi$ to Weyl's character formula and use Lemma~\ref{laDenom}
for $\theta = \rho$ and $\theta = \lambda+\rho$. Since 
$\langle \rho, \alpha^\vee\rangle > 0$ for all $\alpha \in \Phi^+$ we have 
$\psi(A_\rho) \neq 0$, this yields the first equality. The second follows
from $\sum_{\alpha \in \Phi^+} \langle \lambda, \alpha^\vee \rangle =
\langle \lambda, 2\rho^\vee \rangle$.
\end{Prf}

\begin{Def}\label{defMultiset}
For a dominant weight $\lambda$ we define the multiset
\[ M(\lambda) = \mset{\langle \lambda+\rho, \alpha^\vee \rangle \mid 
\alpha \in \Phi^+}. \]
\end{Def}

Since $\lambda + \rho = \sum_i c_i \omega_i$ with all $c_i \geq 1$, 
the multiset $M(\lambda)$ consists of $|\Phi^+|$ positive integers.

\begin{La}\label{laPolys}
For a finite multiset $I$ of positive integers let $P_I = \prod_{i \in I}
(t^i - 1) \in \C[t]$. Then $P_I$ determines $I$.
\end{La}
\begin{Prf}
We use induction on $|I|$, the case $I = \emptyset$ being trivial. 
The zeroes of $t^i - 1$ are the $i$-th roots of unity, and a primitive
$m$-th root of unity is a zero of $t^i-1$ if and only if $m \mid i$. 
Therefore, the largest element $m$ of $I$ is the largest integer such that
a primitive $m$-th root of unity is a zero of $P_I$. Removing one copy of
$m$ from $I$ yields a multiset $I'$ with $P_{I'} = P_I / (t^m-1)$ and we
can apply the induction hypothesis.
\end{Prf}

\begin{Cor}\label{corMultiset}
Let $\lambda$, $\mu$ be dominant weights and $A$ be a principal
$A_1$-subgroup of $G$. Then the restrictions of $V(\lambda)$ and $V(\mu)$
to $A$ are isomorphic if and only if $M(\lambda) = M(\mu)$. 

Equivalently, this holds if and only if 
$\mset{\langle \lambda+\rho, \alpha^\vee \rangle \mid \alpha \in \Phi} = 
\mset{\langle \mu+\rho, \alpha^\vee \rangle \mid \alpha \in \Phi}$.
\end{Cor}
\begin{Prf}
By Section~\ref{secPrincipal} we have to show that $\chi_\lambda(2\rho^\vee(t)) =
\chi_\mu(2\rho^\vee(t))$ if and only if $M(\lambda) = M(\mu)$.

If $M(\lambda) = M(\mu)$ then also 
$\langle \lambda, 2\rho^\vee \rangle = \sum_{m \in M(\lambda)} m - 
\sum_{\alpha \in \Phi^+} \langle \rho, \alpha^\vee \rangle = 
\langle \mu, 2\rho^\vee \rangle$, and Theorem~\ref{thmSerre} shows that the
Laurent polynomials are equal.

Conversely, assume $\chi_\lambda(2\rho^\vee(t)) = \chi_\mu(2\rho^\vee(t))$.
The weights of $V(\lambda)$ are of the form $\lambda - \sum_i n_i \alpha_i$
with all $n_i \in \Z_{\geq 0}$, and $\langle \alpha_i, 2\rho^\vee \rangle =
2 > 0$. So the largest power of $t$ occurring in $\chi_\lambda(2\rho^\vee(t))$
is $t^{\langle \lambda, 2\rho^\vee \rangle}$, and we get $\langle \lambda,
2\rho^\vee \rangle = \langle \mu, 2\rho^\vee \rangle$. Multiplying the
formula in Theorem~\ref{thmSerre} with $t^{\langle \lambda, 2\rho^\vee
\rangle}$ and with the product of the denominators we get
\[ P_{2M(\lambda)} = P_{2M(\mu)}, \]
where $2M(\lambda)$ is the multiset of all $2m$, $m \in M(\lambda)$, and
$P_I$ is defined as in Lemma~\ref{laPolys}. Now Lemma~\ref{laPolys} shows
$M(\lambda) = M(\mu)$.

The second statement follows because 
$\mset{\langle \lambda+\rho, \alpha^\vee \rangle \mid \alpha \in \Phi}$ 
is the union of $M(\lambda)$ and $-M(\lambda)$, and all elements of
$M(\lambda)$ are positive.
\end{Prf}

\begin{Rem}\label{remDegree}
Clearly, if the restrictions of $V(\lambda)$ and $V(\mu)$ to a principal
$A_1$-subgroup are isomorphic then $V(\lambda)$ and $V(\mu)$ have the same
degree. (The degree of $V(\lambda)$ can be computed from $M(\lambda)$ by
Weyl's dimension formula.)

The converse is far from true: in~\cite{luebeck} it is shown that every
simple simply-connected group of rank $> 1$ has infinitely many pairs of
irreducible representations of the same degree which are not related by a
graph automorphism of the group. 
\end{Rem}

\section{Exceptional groups}\label{secResults}

By Corollary~\ref{corMultiset}, Theorem~\ref{thmMain} for exceptional groups
is equivalent to the following theorem.

\begin{Thm}\label{thmMultiset}
Let $G$ be of exceptional type $G_2$, $F_4$, $E_6$, $E_7$ or $E_8$, and let
$\lambda$, $\mu$ be dominant weights with $M(\lambda) = M(\mu)$. Then
$\lambda = \mu$, or $G$ is of type $E_6$ and $\mu$ is the image of
$\lambda$ under the non-trivial graph automorphism.
\end{Thm}

In the following, we write a weight $\lambda + \rho = c_1 \omega_1 + \ldots
+ c_l \omega_l$ with all $c_i \geq 1$, and we write a positive coroot
$\alpha^\vee = b_1 \alpha_1^\vee + \ldots + b_l \alpha_l^\vee$ in short form
as the string of digits $b_1 b_2 \ldots b_l$. Then $\langle \lambda + \rho,
\alpha^\vee \rangle = c_1 b_1 + \ldots + c_l b_l$.

\begin{Def}\label{defDominate}
Let $\alpha^\vee, \beta^\vee \in \Phi^\vee_+$. We say that $\alpha^\vee$
\emph{dominates} $\beta^\vee$ if $\alpha^\vee \neq \beta^\vee$ and 
$\alpha^\vee - \beta^\vee$ is a linear combination of the simple coroots 
with non-negative coefficients.
\end{Def}

Clearly, if $\alpha^\vee$ dominates $\beta^\vee$ then $\langle \lambda +
\rho, \alpha^\vee \rangle > \langle \lambda + \rho, \beta^\vee \rangle$ 
for all dominant weights $\lambda$.

For a subset $S \subseteq \Phi^\vee_+$ and a dominant weight $\lambda$ we
write $M_S(\lambda) = \mset{\langle \lambda+\rho, \beta^\vee \rangle \mid
\beta^\vee \in S}$, so that $M(\lambda) = M_{\Phi^\vee_+}(\lambda)$.

\begin{La}\label{laChain}
Let $S \subseteq \Phi^\vee_+$ with $|S| = k$, let $\lambda$ be a dominant
weight and let $m_1 \geq m_2 \geq \ldots \geq m_k$ be the elements of 
$M_S(\lambda)$. Enumerate $S = \{\beta_1^\vee, \ldots, \beta_k^\vee\}$
such that $m_i = \langle \lambda + \rho, \beta_i^\vee \rangle$ for
$1 \leq i \leq k$ (this enumeration may not be unique). 

Then for every $1 \leq j \leq k$ the coroot
$\beta_j^\vee$ is not dominated by any of $\beta_{j+1}^\vee, \ldots,
\beta_k^\vee$.
\end{La}
\begin{Prf}
If $\beta_i^\vee$ dominates $\beta_j^\vee$ for some 
$i > j$ then $m_i > m_j$ by the remark above, a contradiction.
\end{Prf}

\begin{Rem}\label{remChain}
Let $S$ and $k$ be as in Lemma~\ref{laChain} and let $1 \leq j \leq k$.

Lemma~\ref{laChain} allows us to enumerate, independently of $\lambda$,
all possibilities for the coroots $\beta_1^\vee, \ldots, \beta_j^\vee$
corresponding to the $j$ largest elements of $M_S(\lambda)$: the coroot
$\beta_1^\vee$ must be a maximal element of $S$ with respect to
domination, that is, it is not dominated by any element of $S$; then
$\beta_2^\vee$ must be a maximal element of $S \setminus
\{\beta_1^\vee\}$, and so on. 

If this yields only one possible sequence then 
$m_i = \langle \lambda+\rho, \beta_i^\vee \rangle$ for $1 \leq i \leq j$
and all dominant weights $\lambda$, and we get $j$ linear equations for 
$\lambda + \rho$.
\end{Rem}

\begin{Rem}\label{remSum}
Since $2\rho^\vee$ is the sum of the positive coroots we have 
\[ \sum_{m \in M(\lambda)} m = \langle \lambda + \rho, 2\rho^\vee \rangle,\]
so $M(\lambda)$ always yields this further linear equation for 
$\lambda + \rho$. 

If some linear equations for $\lambda + \rho$ are known then we can compute
$\langle \lambda+\rho, \beta^\vee \rangle$ for all $\beta^\vee$ in the
set $S'$ of positive coroots which lie in the span of the corresponding
elements of $\Q \otimes Y$. Removing these values from $M(\lambda)$ we obtain the
multiset $M_S(\lambda)$ for $S = \Phi^\vee_+ \setminus S'$, to which we
can apply Remark~\ref{remChain} again.
\end{Rem}

\subsection{Types $G_2$ and $F_4$}\label{secG2F4}

Recall that we use the labelling of the simple roots
from~\cite[Planches]{bou}.

\begin{Pro}\label{proG2F4}
Theorem~\ref{thmMultiset} holds for $G$ of type $G_2$ and $F_4$. More
precisely, let $m_1 \geq m_2 \geq \ldots$ be the elements of $M(\lambda)$.
\begin{enumerate}
\item If $G$ is of type $G_2$ then $c_1 = m_1 - m_2$ and 
$c_2 = (2m_2 - m_1)/3$.
\item If $G$ is of type $F_4$ then $c_4 = m_1 - m_2$, $c_3 = m_2 - m_3$, 
$c_2 = m_3 - m_4$ and $c_1 = (m_4 - 3c_2 - 2c_3 - c_4)/2$.
\end{enumerate}
\end{Pro}
\begin{Prf}
The positive coroots of $G_2$ are
\[ 10,\ 01,\ 11,\ 12,\ 13,\ 23. \]
Remark~\ref{remChain} applied to $\Phi^\vee_+$ and $j = 2$ yields only one
possible sequence, $23$, $13$: the coroot $23$ dominates all other positive
coroots, and $13$ dominates $10$, $01$, $11$ and $12$. 
So $m_1 = 2c_1 + 3c_2$ and $m_2 = c_1 + 3c_2$, and (a) follows.

The positive coroots of $F_4$ are
\[ \begin{array}{llllllll}
1000, & 0100, & 0010, & 0001, & 1100, & 0110, & 0011, & 1110, \\
0210, & 0111, & 1210, & 1111, & 0211, & 2210, & 1211, & 0221, \\
2211, & 1221, & 2221, & 1321, & 2321, & 2421, & 2431, & 2432.
\end{array} \]
Remark~\ref{remChain} applied to $\Phi^\vee_+$ and $j = 4$ yields only one
possible sequence, $2432$, $2431$, $2421$, $2321$: each of these coroots
dominates all positive coroots which appear before it in the list above.
So we get
$m_1 = 2c_1 + 4c_2 + 3c_3 + 2c_4$,
$m_2 = 2c_1 + 4c_2 + 3c_3 + c_4$,
$m_3 = 2c_1 + 4c_2 + 2c_3 + c_4$ and
$m_4 = 2c_1 + 3c_2 + 2c_3 + c_4$, and (b) follows.
\end{Prf}

\subsection{Type $E_8$}\label{secE8}

\begin{Pro}\label{proE8}
Theorem~\ref{thmMultiset} holds for $G$ of type $E_8$.
\end{Pro}
\begin{Prf}
There are $120$ positive coroots 
and
\[ 2\rho^\vee = (92, 136, 182, 270, 220, 168, 114, 58). \]

Let $m_1 \geq m_2 \geq \ldots$ be the elements of $M(\lambda)$.
Remark~\ref{remChain} applied to $\Phi^\vee_+$ and $j = 6$ yields only one
possible sequence
\[ 23465432,\ 23465431,\ 23465421,\ 23465321,\ 23464321,\ 23454321. \]
This shows $c_8 = m_1 - m_2$, $c_7 = m_2 - m_3$, $c_6 = m_3 - m_4$, 
$c_5 = m_4 - m_5$, $c_4 = m_5 - m_6$, and 
$m_6 = 2c_1 + 3c_2 + 4c_3 + 5c_4 + 4c_5 + 3c_6 + 2c_7 + c_8$, so we know
$2c_1 + 3c_2 + 4c_3$. By Remark~\ref{remSum} we also know 
$\langle \lambda+\rho, 2\rho^\vee\rangle$ and so 
$92c_1 + 136c_2 + 182c_3$.

The six coroots above and $2\rho^\vee$ are linearly independent, and
exactly $42$ positive coroots lie in their span. As explained in
Remark~\ref{remSum}, we can compute the values in $M(\lambda)$ for these
$42$ coroots, and we obtain $M_S(\lambda)$ for the set $S$ of the remaining
$78$ positive coroots. Let $m'_1 \geq m'_2 \geq \ldots$ be the elements of
$M_S(\lambda)$. Remark~\ref{remChain} applied to $S$ and $j=4$ yields three
possible sequences: the first three coroots are $13354321$, $23354321$,
$22454321$ in some order, and the fourth coroot is always $12354321$. This
shows $m'_4 = c_1 + 2c_2 + 3c_3 + 5c_4 + 4c_5 + 3c_6 + 2c_7 + c_8$, so we
know $c_1 + 2c_2 + 3c_3$.

Since 
\[ \det \begin{pmatrix} 2 & 3 & 4 \\ 92 & 136 & 182 \\ 1 & 2 & 3
\end{pmatrix} = -2 \neq 0, \]
we can also compute $c_1$, $c_2$ and $c_3$. So $M(\lambda)$ determines 
$\lambda + \rho$ and hence $\lambda$.
\end{Prf}

\subsection{Type $E_6$}\label{secE6}

Let $\sigma$ be the non-trivial graph automorphism of type $E_6$. It acts on
$X$ and $Y$ by $\omega_1 \leftrightarrow \omega_6$, $\omega_3
\leftrightarrow \omega_5$, $\alpha_1^\vee \leftrightarrow \alpha_6^\vee$, 
$\alpha_3^\vee \leftrightarrow \alpha_5^\vee$, and it fixes $\omega_2$,
$\omega_4$, $\alpha_2^\vee$ and $\alpha_4^\vee$. In particular, $\sigma$
permutes $\Phi^\vee_+$, fixes $\rho$, and $\langle \sigma(x), \sigma(y)
\rangle = \langle x, y \rangle$ for all $x \in X$, $y \in Y$. Hence 
$M(\sigma(\lambda)) = M(\lambda)$ for all dominant weights $\lambda$, and
the exception in Theorem~\ref{thmMultiset} really occurs.

\begin{Pro}\label{proE6}
Theorem~\ref{thmMultiset} holds for $G$ of type $E_6$.
\end{Pro}
\begin{Prf}
There are $36$ positive coroots and 
$2\rho^\vee = (16, 22, 30, 42, 30, 16)$.

Let $\lambda$, $\mu$ be dominant weights with $M(\lambda) = M(\mu)$. Since
$M(\sigma(\lambda)) = M(\lambda)$ and $\sigma$ interchanges $c_3$ and $c_5$,
we may replace $\lambda$ and $\mu$ by their images under $\sigma$, if
necessary, and assume that $c_3 \leq c_5$ holds for both of them. We show
that $M(\lambda)$ determines $c_2$, $c_3$, $c_4$, $c_5$ and $c_1 + c_6$, and
the set $\{c_1, c_6\}$, and that $\mu = \lambda$ or $\mu = \sigma(\lambda)$.

Let $m_1 \geq m_2 \geq \ldots$ be the elements of $M(\lambda)$.
Remark~\ref{remChain} applied to $\Phi^\vee_+$ and $j=4$ yields two possible
sequences: the first three coroots are $122321$, $112321$, $112221$, and 
the fourth one is $111221$ or $112211$. So $m_4$ is the larger of the
values of $\lambda + \rho$ on $111221$ and $112211$. Since $111221 - 112211
= -\alpha_3^\vee + \alpha_5^\vee$ and $c_3 \leq c_5$ this is the value on
$111221$. We get $c_2 = m_1 - m_2$, $c_4 = m_2 - m_3$, $c_3 = m_3 - m_4$
and $m_4 = c_1 + c_2 + c_3 + 2c_4 + 2c_5 + c_6$. Together with 
$\langle \lambda+\rho, 2\rho^\vee\rangle = 16(c_1+c_6) + 22c_2 + 30(c_3 +
c_5) + 42c_4$ from Remark~\ref{remSum} this determines $c_5$ and 
$c_1 + c_6$.

With this information we can compute the values in $M(\lambda)$ for the
$20$ positive coroots with $b_1 = b_6$.
Let $S$ be the set of the
remaining $16$ positive coroots:
\[ \begin{array}{cccccccc}
000001 & 100000 & 000011 & 101000 & 000111 & 101100 & 001111 & 010111 \\
101110 & 111100 & 011111 & 111110 & 011211 & 111210 & 011221 & 112210
\end{array} \]
By Remark~\ref{remSum} we know $M_S(\lambda)$. Each element of $S$ is
equal to or dominates $100000$ or $000001$. So the smallest element $m$
of $M_S(\lambda)$ is $\min\{c_1, c_6\}$, and we get 
$\{c_1, c_6\} = \{m, (c_1+c_6) - m\}$.

This shows that $\mu = \lambda$, or $c_1 \neq c_6$ and $\mu$ is obtained
from $\lambda$ by interchanging $c_1$ and $c_6$. In the second case we may
assume that $c_1 < c_6$. The maximal elements of $S$ with respect to
domination are $011221$ and $112210$, so the largest element of
$M_S(\lambda) = M_S(\mu)$ is the largest of the values on these two coroots. 
For $\lambda$ we have 
\[ \langle \lambda+\rho, 011221 \rangle - 
   \langle \lambda+\rho, 112210 \rangle = (c_6 - c_1) + (c_5 - c_3) > 0, \]
so the largest element of $M_S(\lambda)$ is $c_2 + c_3 + 2c_4 + 2c_5 +
c_6$. For $\mu$ the value on $011221$ is $c_2 + c_3 + 2c_4 + 2c_5 + c_1$,
which is smaller. Therefore, the largest element of $M_S(\mu)$ is the value
of $\mu + \rho$ on $112210$, that is $c_6 + c_2 + 2c_3 + 2c_4 + c_5$.
Comparing the two expressions yields $c_3 = c_5$, and so $\mu =
\sigma(\lambda)$.
\end{Prf}

\subsection{Type $E_7$}\label{secE7}

In this case the method of this section only gives partial results. We
show how to compute from $M(\lambda)$ a set of at most three dominant
weights which contains $\lambda$. A full proof of this case is given in
Theorem~\ref{thmInv}.

There are $63$  positive coroots 
and $2\rho^\vee = (34, 49, 66, 96, 75, 52, 27)$.

Remark~\ref{remChain} applied to $\Phi^\vee_+$ and $j = 6$ yields
three possible sequences: the first four coroots are $2234321$, $1234321$,
$1224321$, $1223321$, and the next two can be
\begin{enumerate}
\item[(1)] $1223221$, $1223211$, which is possible if and only if 
$c_2 \geq c_5 + c_6$,
\item[(2)] $1223221$, $1123321$, which is possible if and only if 
$c_5 \leq c_2 \leq c_5 + c_6$, or
\item[(3)] $1123321$, $1223221$, which is possible if and only if 
$c_2 \leq c_5$.
\end{enumerate}
In each case these six coroots and $2\rho^\vee$ are linearly independent.
So, by Remark~\ref{remSum}, there are at most three dominant weights with a
given multiset $M(\lambda)$, and they can easily be computed from the sum
and the six largest elements of $M(\lambda)$ (often one or two of the
solutions have non-integer or negative entries or do not fulfill the
mentioned conditions).

In general, it seems difficult to decide by considering largest elements
in $M(\lambda)$ which of these
candidates is the correct one. For example, for $\lambda = \omega_1 +
3\omega_5 + 3\omega_6 + \omega_7$ (case (3)) and $\mu = \omega_1 +
4\omega_2 + 2\omega_6 + 4\omega_7$ (case (1)) the sums of the elements of
$M(\lambda)$ and $M(\mu)$ are equal (to $841$) and also the $20$ largest
elements of $M(\lambda)$ and $M(\mu)$ coincide, they are
\[ 35, 33, 32, 31, 30, 27, 26, 25, 24, 23, 22, 22, 21, 21, 20, 20, 19, 18,
18, 17. \]
The $21$st largest elements are $16$ and $17$, respectively.

\section{Classical types}\label{secClassical}

The elementary methods of Section~\ref{secResults} also work for the
simple groups of classical type of small rank. We keep the notation
$\lambda + \rho = \sum_i c_i \omega_i$ and write coroots as strings of
digits. The non-trivial graph automorphism of type $A_l$ reverses the
sequence $(c_1, \ldots, c_l)$, and the graph automorphisms of type $D_4$
permute $c_1$, $c_3$, $c_4$ arbitrarily.

We now prove the statement in~\cite[\S10]{npp} that for all simple groups of
rank at most $4$ the restriction to a principal $A_1$-subgroup determines
the irreducible representation up to graph automorphisms. Types $G_2$ and
$F_4$ were done in Proposition~\ref{proG2F4}.
\begin{Thm}\label{thmClassical}
Let $G$ be of type $A_l$ with $2 \leq l \leq 4$, $B_l$ or $C_l$ with 
$2 \leq l \leq 4$, or $D_4$, and let $\lambda$, $\mu$ be dominant weights
with $M(\lambda) = M(\mu)$. Then $\mu$ is the image of $\lambda$ under a
graph automorphism. In particular, $\lambda = \mu$ for the types $B_l$ and
$C_l$.
\end{Thm}
\begin{Prf}
Let $m_1 \geq m_2 \geq \ldots$ be the elements of $M(\lambda)$. In each
case we only sketch the argument, the lists of positive coroots are easily
written down.

\emph{Types $B_2$, $C_2$, $B_3$, $C_3$.} Remark~\ref{remChain} applied to
$\Phi^\vee_+$ and $j=2$ yields only one possible sequence, namely $21$,
$11$ for $B_2$, $12$, $11$ for $C_2$, $221$, $121$ for $B_3$ and $122$, $112$
for $C_3$. For $B_2$ and $C_2$ this yields two independent equations, for
$B_3$ and $C_3$ we get a third one from Remark~\ref{remSum}, since 
$2\rho^\vee$ has coefficients $(6,10,6)$ and $(5,8,9)$, respectively.

\emph{Type $B_4$.} Here $2\rho^\vee$ has coefficients $(8,14,18,10)$. The only possible
sequence for $j = 2$ is $2221$, $1221$, which yields $c_1 = m_1 - m_2$ and 
$c_1 + 2c_2 + 2c_3 + c_4$. Together with Remark~\ref{remSum} we know the
values on the positive coroots $1000$, $0221$, $1221$, $2221$. For the set
$S$ of the remaining $12$ positive coroots, $1121$ dominates all other
elements of $S$, so the largest element of $M_S(\lambda)$ is $c_1 + c_2 +
2c_3 + c_4$. This gives a fourth independent linear equation.

\emph{Type $C_4$.} Here $2\rho^\vee$ has coefficients $(7,12,15,16)$. 
The largest two entries of $M(\lambda)$ correspond to $1222$, $1122$ and the
third to $1112$ or $0122$. So, if $\lambda \neq \mu$, then after possibly
interchanging $\lambda$ and $\mu$ we have
$\lambda+\rho = c_1\omega_1 + c_2\omega_2 + c_3\omega_3 + c_4\omega_4$
and $\mu + \rho = c_3\omega_1 + c_2\omega_2 + c_1\omega_3 +
(c_4+(c_3-c_1)/2)\omega_4$.

Comparing $M(\lambda)$ and $M(\mu)$, $12$ of the values coincide. With $d =
(c_3 - c_1)/2$ the remaining values are
\[ c_4 < c_3 + c_4 < c_2 + c_3 + c_4 < c_1 + c_2 + c_3 + c_4 \]
for $\lambda + \rho$ and
\[ c_4 + d < c_1 + c_4 + d < c_1 + c_2 + c_4 + d < c_1 + c_2 + c_3 + c_4 + d \]
for $\mu + \rho$. The smallest values are equal if and only if $d = 0$,
that is $c_1 = c_3$, and so $\lambda = \mu$.

\emph{Type $D_4$.} Here $2\rho^\vee$ has coefficients $(6,10,6,6)$. Applying a graph
automorphism we may assume $c_1 \leq c_3 \leq c_4$, and similarly for
$\mu$. The only possible sequence for $j=2$ is $1211$, $1111$, so we know
$c_2 = m_1 - m_2$ and $c_1 + c_3 + c_4$. The next largest value is on one
of $0111$, $1110$, $1101$, and by our assumption it is $c_2 + c_3 + c_4$;
this yields $c_1$. Now we know the values on the six positive coroots with
$b_3 = b_4$. Each of the remaining six positive coroots is equal to or
dominates $0010$ or $0001$, so the smallest of their values is $c_3$, and
then we also know $c_4$.

\emph{Types $A_2$, $A_3$, $A_4$.} Applying the graph automorphism we may
assume $c_1 \leq c_l$, and similarly for $\mu$. For $A_2$ the multiset
$M(\lambda) = \mset{c_1, c_2, c_1 + c_2}$ determines $c_1 \leq c_2$. For
$A_3$ and $A_4$ the largest element $m_1$ is the sum of all $c_i$ and
$m_2 = m_1 - c_1$, this yields $c_1$. By Remark~\ref{remSum} we also know
$\sum_{m \in M(\lambda)} m$, which is $3m_1 + c_2$ for $A_3$ and $4m_1 +
2(c_2 + c_3)$ for $A_4$. This determines $\lambda$ for $A_3$.

For $A_4$ we now know $c_1$, $c_4$ and $s = c_2 + c_3$, and so the values on
the six positive coroots with $b_2 = b_3$. The remaining values are
$\mset{c_2, c_3, c_1 + c_2, c_3 + c_4}$ on the coroots $0100$, $0010$,
$1100$, $0011$. Each of these is equal to or dominates $0100$ or $0010$, so
the smallest of these values is $m = \min\{c_2, c_3\}$, and $\{c_2, c_3\}
= \{m, s - m\}$. If $c_2 \neq c_3$ and both choices yield the same
multiset, then $\mset{c_1 + m, c_4 + s - m} = \mset{c_1 + s - m, c_4 +
m}$. Since $m < s - m$ this implies $c_1 = c_4$, and the two choices are
related by the graph automorphism.
\end{Prf}

\begin{Rem}\label{remAl}
In type $A_l$ the values $\langle \lambda+\rho, \alpha^\vee \rangle$ are
the sums $c_i + c_{i+1} + \ldots + c_j$, $1 \leq i \leq j \leq l$. So, if we
set $p_0 = 0$ and $p_k = c_1 + \ldots + c_k$ for $1 \leq k \leq l$, then
$M(\lambda)$ is the multiset of distances $p_j - p_i$, $0 \leq i < j \leq
l$, between the $l+1$ points $p_0 < p_1 < \ldots < p_l$, and the graph
automorphism corresponds to the reflection $p \mapsto p_l - p$. So the
question whether $M(\lambda)$ determines $\lambda$ up to graph automorphism
translates to the classical question on \emph{homometric sets of integers}, 
that is, sets with the same multiset of distances.

For $l = 5$ there are well-known examples of homometric sets which are not
related by shift or reflection. For example, for positive integers $a$ and
$c$ set $\lambda+\rho = (a, c, 3a+c, 2a, 2a+c)$ and 
$\mu+\rho = (a, 4a+c, c, 2a, a+c)$. Then it is easily checked that
$M(\lambda) = M(\mu)$. For $a=1$, $c=3$ this is the example of 
Bloom~\cite{bloom}. 

Using certain polynomial identities as in~\cite{rosseym} one can also
construct examples for $A_l$ with $l \geq 6$. For example for
$\lambda+\rho = (4, 1, 2l - 5, 3, 1, \ldots, 1, 2, 1)$ and
$\mu+\rho = (1, 4, 2l - 5, 1, 2, 1, \dots, 1, 3)$ 
(the $l-6$ entries between $3$ and $2$, or between $2$ and $3$ respectively,
are all $1$) we have $M(\lambda) = M(\mu)$.
\end{Rem}

%
%

\section{Invariant theory}\label{secInvariants}

In this section we give a proof of Theorem~\ref{thmMain} for several types,
including $E_7$, using the invariant theory of the Weyl group. 
This approach was also used by Serre in~\cite{serreA1} to show
Theorem~\ref{thmMain} for types $E_7$ and $E_8$.

Let $V = \C \otimes X$ and let $R = \C[x_1, \ldots, x_l]$ be the ring of
polynomial functions on $V$, where $x = \sum_i x_i \omega_i \in V$. The
Weyl group $W$ acts on $V$, and on $R$ by $(w f)(x) = f(w^{-1}x)$. By
Chevalley's theorem the ring of invariants $R^W$ is generated by $l$
algebraically independent homogeneous polynomials whose degrees $d_1 \leq
\ldots \leq d_l$ are uniquely determined by $W$~\cite[Ch.~3]{humrg}; these
are the degrees of $W$, see~\cite[3.7]{humrg} for a table.

\begin{La}\label{laBasic}
Let $f_1, \ldots, f_l \in R^W$ be algebraically independent homogeneous
polynomials of degrees $d_1, \ldots, d_l$. Then $R^W = \C[f_1, \ldots,
f_l]$.
\end{La}
\begin{Prf}
The subring $\C[f_1, \ldots, f_l] \subseteq R^W$ and $R^W$ are graded
by degree, and both are polynomial rings in generators of degrees $d_1,
\ldots, d_l$. So their homogeneous components of each degree have the same
finite dimension, and the inclusion is an equality.
\end{Prf}

For a positive integer $d$ we define
\[ S_d(x) = \sum_{\alpha \in \Phi} \langle x, \alpha^\vee \rangle^d \in R. \]
Since $\langle w^{-1}x, \alpha^\vee \rangle = \langle x, w(\alpha^\vee)
\rangle$ and $W$ permutes $\Phi^\vee$, the polynomial $S_d$ is $W$-invariant
and homogeneous of degree $d$. For odd $d$ the summands for $\alpha$ and
$-\alpha$ cancel, so $S_d = 0$. For even $d$ and a dominant weight
$\lambda$ we have
\[ S_d(\lambda + \rho) = 2 \sum_{m \in M(\lambda)} m^d, \]
so $M(\lambda)$ determines the values of all $S_d$ at $\lambda + \rho$.

To check that $S_{d_1}, \ldots, S_{d_l}$ are algebraically independent we
use the Jacobian criterion, see~\cite[3.10]{humrg}: this holds if and only
if the Jacobian determinant $\det(\partial S_{d_i} / \partial x_j)$ is not
the zero polynomial, and for this it is enough to find one point where it
does not vanish. We evaluate at $x = \rho$. For $h \geq 1$ let $n_h$ be the
number of positive coroots of height $h$, and let
\[ \sigma_h = \sum_{\alpha^\vee \in \Phi^\vee_+,\ \hgt(\alpha^\vee) = h}
\alpha^\vee = \sum_{j=1}^l s_{hj}\, \alpha_j^\vee \in Y \]
be the sum of these coroots. Let $H$ be the height of the highest coroot.
Since $\langle \rho, \alpha^\vee \rangle = \hgt(\alpha^\vee)$ for
$\alpha^\vee \in \Phi^\vee_+$, we get for even $d$
\[ S_d(\rho) = 2 \sum_{h=1}^H n_h\, h^d > 0, \]
in particular $S_d \neq 0$. Furthermore, $\partial \langle x, \alpha^\vee
\rangle^d / \partial x_j = d\, b_j \langle x, \alpha^\vee \rangle^{d-1}$
for $\alpha^\vee = \sum_j b_j \alpha_j^\vee$, and since for even $d$ the
summands for $\alpha$ and $-\alpha$ coincide we get
\[ \frac{\partial S_d}{\partial x_j}(\rho) = 2d \sum_{h=1}^H h^{d-1}\,
s_{hj}. \]
So, if all $d_i$ are even, the Jacobian matrix at $\rho$ is the product of
the diagonal matrix with entries $2d_1, \ldots, 2d_l$ and the integer
matrix
\[ A = P \cdot \Sigma \quad \textrm{with} \quad P = \left(h^{d_i - 1}\right)_{1
\leq i \leq l,\ 1 \leq h \leq H} \quad \textrm{and} \quad 
\Sigma = \left(s_{hj}\right)_{1 \leq h \leq H,\ 1 \leq j \leq l}. \]
To show $\det(A) \neq 0$ it suffices to find a prime $p$ such that
$\det(A) \not\equiv 0 \pmod p$; this is a fast computation even for large
$l$.

For example, for type $E_7$ we have $H = 17$, the degrees are $2$, $6$, $8$, $10$,
$12$, $14$, $18$, and we find $\det(A)
\equiv 5 \pmod{19}$ (we do not need to compute  $\det(A) \approx 1.025 \cdot
10^{64}$ or even the Jacobian determinant as a huge polynomial). 
For the types $G_2$, $F_4$ and $E_8$, with degrees
$2, 6$, and $2, 6, 8, 12$, and $2, 8, 12, 14, 18, 20, 24, 30$, the
smallest primes $p$ with $\det(A) \not\equiv 0 \pmod p$ are $7$, $13$ and
$31$, respectively.

To deal with the classical types we do not use the Jacobian criterion
directly. Instead we show that another set of well-known generating
invariants (for which the Jacobian criterion was used) 
can be expressed as polynomials in the $S_{d_i}$. 

For this we first recall the description in~\cite[Planches]{bou} of 
the coroot system and the
corresponding action of the Weyl group in types $B_l$, $C_l$ and $D_l$.
There is a basis $\varepsilon_1, \ldots, \varepsilon_l$ of 
$\Q \otimes Y$ such that the positive coroots are
\[ \begin{array}{ll}
\varepsilon_i \pm \varepsilon_j & \textrm{for }1 \leq i < j \leq l,\\
2\varepsilon_i &  \textrm{for }1 \leq i \leq l \textrm{ (only in type
$B_l$)},\\
\varepsilon_i &  \textrm{for }1 \leq i \leq l \textrm{ (only in type
$C_l$)}.
\end{array}
\]

Let $y_i \in R$ be the linear form $y_i(x) = \langle
x, \varepsilon_i \rangle$.
Let $W'$ be the group of linear maps on $\Q \otimes Y$ which permute the
$\{\pm \varepsilon_i\mid\; 1\leq i \leq l\}$. 
For the types $B_l$ and
$C_l$ this is exactly the Weyl group  $W = W'$. 
For type $D_l$ the Weyl group $W$ is the subgroup of
index $2$ of $W'$ consisting of the elements with an even number of sign
changes, and $W' = W \cup W\tau$, where $\tau$ just changes the sign of
$\varepsilon_l$. 

In~\cite[3.12]{humrg} it is shown that the ring of invariants $R^{W'}
= \C[q_1, \ldots, q_l]$ is generated by  $q_a = \sum_{i=1}^l y_i^{2a}$
for $1\leq a \leq l$. We also write $q_0 = l$.

\begin{La}\label{laBC}
Let $G$ be of type $B_l$ or $C_l$ with $l \geq 2$, or of type $D_l$ with $l
\geq 4$. 
\begin{enumerate}
\item For $n \geq 1$ we have
\[ S_{2n} = \gamma_n\, q_n + 2 \sum_{k=1}^{n-1} \binom{2n}{2k} q_k
q_{n-k}, \]
where
\[ \gamma_n = \begin{cases} 4l + 4^n & \textrm{for type } B_l, \\ 
4l + 2 - 4^n & \textrm{for type } C_l, \\
4l - 4^n & \textrm{for type } D_l. \end{cases} \]
\item For the types $B_l$ and $C_l$, and for type $D_l$ when $l$ is not a
power of $4$, we have $R^{W'} = \C[S_2, S_4, \ldots, S_{2l}]$.
\item In case $D_l$ with $l = 4^m$ the invariant $S_{2m + 2}$ is a
polynomial in the $S_d$ with $d < 2m+2$. So, the $S_d$, $d \geq 1$, do not
generate $R^{W'}$ in this case.
\end{enumerate}
\end{La}
\begin{Prf}
Let $d = 2n$ with $n \geq 1$. We compute $S_d$ in terms of the $q_a$. 
Recall that $\langle x, \varepsilon_i \pm \varepsilon_j \rangle = (y_i \pm
y_j)(x)$. We have
\[ (y_i + y_j)^d + (y_i - y_j)^d = 2 \sum_{k=0}^n \binom{2n}{2k}
y_i^{2n-2k} y_j^{2k}. \]
Expanding the products $q_{n-k} q_k$ 
and using $\binom{2n}{2k} = \binom{2n}{2n-2k}$ we get
\[ \sum_{1 \leq i < j \leq l} \left((y_i + y_j)^d + (y_i - y_j)^d\right) =
\sum_{k=0}^n \binom{2n}{2k} \left(q_{n-k} q_k - q_n\right). \]
The negative coroots give the same contribution. The summands for $k = 0$
and $k = n$ (with $q_0 = l$),
$\sum_{k=0}^n \binom{2n}{2k} = 2^{2n-1}$ and $S_d = 2\sum_{\alpha^\vee
\in \Phi^\vee_+} \langle x, \alpha^\vee \rangle^d$ yield~(a) in case $D_l$.
In types $B_l$ and $C_l$ the
coroots $\pm 2\varepsilon_i$, respectively $\pm \varepsilon_i$, contribute
$2 \cdot 4^n q_n$ in type $B_l$, respectively $2q_n$ in type $C_l$, to
$S_{d}$. This concludes the proof of~(a).

In type $B_l$ we clearly have $\gamma_n \neq 0$, in type $C_l$ we have
$\gamma_n \neq 0$ because $4l + 2 \not\equiv 0 \pmod 4$, and in type $D_l$
we have $\gamma_n \neq 0$ if and only if $l \neq 4^{n-1}$. So in the cases of (b)
all $\gamma_n$ are non-zero, and by induction on $n$ it follows from (a)
that $\C[S_2, \ldots, S_{2n}] = \C[q_1, \ldots, q_n]$ for $1 \leq n \leq l$.
The case $n = l$ is our claim.

In case $D_l$ with $l = 4^m$ and $n = m+1$ we see from part~(a) that
$\gamma_n = 0$. Here $\gamma_k \neq 0$ for $k \leq m$, so as in the proof
of~(b) we have $\C[S_2, \ldots, S_{2m}] = \C[q_1, \ldots, q_m]$, and by~(a)
the invariant $S_{2n}$ is a polynomial in $q_1, \ldots, q_m$, hence in the
$S_d$ with $d < 2n$. Since $S_d = 0$ for odd $d$, the homogeneous component
of degree $2n$ of $\C[S_d \mid d \geq 1]$ is contained in $\C[q_1, \ldots,
q_m]$, so it does not contain $q_n$ (note that $n \leq l$).
\end{Prf}

\begin{Thm}\label{thmInv}
Let $G$ be of type $G_2$, $F_4$, $E_7$, $E_8$, or of type $B_l$ or $C_l$
with $l \geq 2$, and let $\lambda$, $\mu$ be dominant weights with
$M(\lambda) = M(\mu)$. Then $\lambda = \mu$.
\end{Thm}
\begin{Prf}
In all these cases the degrees $d_1, \ldots, d_l$ of $W$ are even, and
the polynomials $S_{d_1}, \ldots, S_{d_l}$ generate $R^W$. For the types
$G_2$, $F_4$, $E_7$, $E_8$ they are algebraically independent by the
computations described above, so this follows from Lemma~\ref{laBasic}.
For the types $B_l$ and $C_l$ this is Lemma~\ref{laBC}(b).
Since $S_d(\lambda + \rho) = S_d(\mu + \rho)$ for all $d$ we get
$f(\lambda+\rho) = f(\mu+\rho)$ for all $f \in R^W$. 

Now we use that the invariants of a finite group separate the orbits under
this group, see~\cite[Cor.~2.3.6]{derksenkemper}. 
So $\lambda+\rho$ and
$\mu+\rho$ are in the same $W$-orbit. Both are dominant and each $W$-orbit
contains a unique dominant weight, see~\cite[1.12]{humrg}.
This shows $\lambda + \rho = \mu + \rho$.
\end{Prf}

In particular, this completes the proof of Theorem~\ref{thmMultiset} for
type $E_7$, and it gives alternative proofs for the types $G_2$, $F_4$, $E_8$
and $B_l$, $C_l$ with $l \leq 4$.

\begin{Rem}\label{remInvLimits}
This argument for Theorem~\ref{thmInv} cannot work when $G$ 
has a non-trivial graph automorphism
$\sigma$: then $\sigma$ permutes $\Phi^\vee$, so $S_d(\sigma(x)) = S_d(x)$
for all $d$, and the $S_d$ do not distinguish $\lambda + \rho$ from the
dominant element $\sigma(\lambda + \rho) = \sigma(\lambda) + \rho$.
Accordingly, the Weyl groups of the types $A_l$ ($l \geq 2$), $D_l$ ($l$
odd) and $E_6$ have an odd degree ($3$, $l$, and $5$, $9$, respectively),
for which $S_d = 0$. For type $D_l$ with $l$ even all degrees $2, 4,
\ldots, 2l-2$ and $l$ are even, but the degree $l$ occurs twice, so there
are two algebraically independent homogeneous invariants of degree $l$,
whereas $S_l$ is only one invariant of this degree.
\end{Rem}

For type $D_l$ we can still show that $M(\lambda)$ determines $\lambda$ up
to the graph automorphism, provided that $l$ is not a power of $4$.

\begin{Pro}\label{proDl}
Let $G$ be of type $D_l$ with $l \geq 5$ not a power of $4$, let $\sigma$
be the non-trivial graph automorphism of $G$, and let $\lambda$, $\mu$ be
dominant weights with $M(\lambda) = M(\mu)$. Then $\mu = \lambda$ or $\mu =
\sigma(\lambda)$.
\end{Pro}
\begin{Prf}
By Lemma~\ref{laBC}(b) the polynomials $S_2, \ldots, S_{2l}$ generate
$R^{W'}$. As in the proof of Theorem~\ref{thmInv} it follows that $\lambda
+ \rho$ and $\mu + \rho$ are in the same $W'$-orbit. The element $\tau$
interchanges the simple coroots $\alpha_{l-1}^\vee = \varepsilon_{l-1} -
\varepsilon_l$ and $\alpha_l^\vee = \varepsilon_{l-1} + \varepsilon_l$ and
fixes $\alpha_i^\vee = \varepsilon_i - \varepsilon_{i+1}$ for $i < l-1$, so
it acts on $X$ as the graph automorphism $\sigma$; in particular it fixes
$\rho$. Since $W' = W \cup W\tau$, the element $\mu + \rho$ is in the
$W$-orbit of $\lambda + \rho$ or of $\sigma(\lambda + \rho) = \sigma(\lambda)
+ \rho$. All these elements are dominant, so $\mu + \rho = \lambda + \rho$
or $\mu + \rho = \sigma(\lambda) + \rho$.
\end{Prf}

\begin{Rem}\label{remD4}
For type $D_4$ the polynomials $S_d$ cannot separate the $W'$-orbits,
because they are also invariant under the graph automorphisms of order $3$,
which are not induced by elements of $W'$. Here $\gamma_2 = 0$ and $S_4$ is
a multiple of $S_2^2$. The result for $D_4$ was shown in
Theorem~\ref{thmClassical}.
\end{Rem}

Finally, we reformulate our problem in case $D_l$ as a special case of a
combinatorial question studied by Moser, Selfridge and Straus, and others. This provides in cases $D_l$, $l = 4^m$,
explicit $\lambda$ and $\mu$ which are not related by a graph automorphism
and for which $M(\lambda) = M(\mu)$.

\subsection*{Relation to Moser's problem}

For a finite multiset $Z$ of numbers let $P(Z)$ be the multiset of the
sums $z + z'$ over all pairs of different positions in $Z$. The question
whether $P(Z)$ determines $Z$ goes back to a problem of
Moser~\cite{moser}. Selfridge and Straus~\cite[Thm.~1, Thm.~2]{selfstraus}
showed that this holds if $|Z|$ is not a power of $2$, and that for $|Z| =
2^{k-1}$ there are counterexamples: for numbers $a_1, \ldots, a_k$ the
multisets of the sums of an even, respectively odd, number of the $a_j$
have the same pairwise sums, see~\cite[\S6]{gordonfs} for this explicit
description.

Our question for type $D_l$ is equivalent to Moser's problem for
symmetric multisets. For a dominant weight $\lambda$ let $y_i = \langle
\lambda + \rho, \varepsilon_i \rangle$, so that $y_1 > \ldots > y_{l-1} >
|y_l|$, and let $A(\lambda) = \mset{\pm y_1, \ldots, \pm y_l}$. Since the
positive coroots are the $\varepsilon_i \pm \varepsilon_j$ with $i < j$,
the pairs of positions with entries $\pm y_i$, $\pm y_j$, $i < j$,
contribute the sums $\pm(y_i + y_j)$, $\pm(y_i - y_j)$, and the pairs
$y_i$, $-y_i$ contribute $0$. So
\[ P(A(\lambda)) = M(\lambda) \cup -M(\lambda) \cup \mset{0, \ldots, 0}
\quad \textrm{(with $l$ zeroes)}, \]
and since all elements of $M(\lambda)$ are positive we get $M(\lambda) =
M(\mu)$ if and only if $P(A(\lambda)) = P(A(\mu))$. On the other hand,
$A(\lambda) = A(\mu)$ if and only if the sets $\{|y_i|\}$ for $\lambda$ and
$\mu$ coincide, that is, if and only if $\mu = \lambda$ or $\mu =
\sigma(\lambda)$, where $\sigma$ is the graph automorphism induced by
$\tau$.

For symmetric multisets all odd power sums vanish, and for $Z = A(\lambda)$
the proof of~\cite[Thm.~1]{selfstraus} becomes the computation in
Lemma~\ref{laBC}(a). This explains why we get a positive result for all $l$
which are not a power of $4$, and not only for $l$ not a power of $2$, as
the result of~\cite{selfstraus} for $|Z| = 2l$ would give. Conversely, the
counterexamples above consist of symmetric multisets if $k$ is even, and
this yields the following proposition.

\begin{Pro}\label{proDlCounter}
Let $m \geq 1$, $k = 2m + 2$ and $l = 4^m = 2^{k-2}$, and let $G$ be of type
$D_l$. Let $E$, respectively $O$, be the set of numbers $s - (2^k-1)/2$,
where $s$ runs through the integers $0 \leq s < 2^k$ whose binary expansion
contains an even, respectively odd, number of digits $1$. Then 
$E = \{\pm y_1, \ldots, \pm y_l\}$ and $O = \{\pm z_1, \ldots, \pm z_l\}$
with $y_1 > \ldots > y_l > 0$ and $z_1 > \ldots > z_l > 0$. Let $\lambda$
and $\mu$ be the weights with $\langle \lambda + \rho, \varepsilon_i \rangle
= y_i$ and $\langle \mu + \rho, \varepsilon_i \rangle = z_i$ for $1 \leq i
\leq l$. Then $\lambda$ and $\mu$ are dominant and $M(\lambda) = M(\mu)$.
For $m \geq 2$ the weight $\mu$ is not the image of $\lambda$ under a graph
automorphism.
\end{Pro}
\begin{Prf}
Since $k$ is even, the map $s \mapsto 2^k - 1 - s$, which interchanges all
binary digits of $s$, preserves the parity of the number of digits $1$. So
$E = -E$ and $O = -O$. The elements of $E$ and $O$ are pairwise distinct
half-integers, in particular non-zero, and $|E| = |O| = 2^{k-1} = 2l$. This
shows the claimed form of $E$ and $O$. For $\lambda + \rho = \sum_i c_i
\omega_i$ we have $c_i = y_i - y_{i+1}$ for $1 \leq i < l$ and $c_l =
y_{l-1} + y_l$. These are positive integers, so $\lambda$ is dominant, and
similarly for $\mu$.

Up to the common shift by $-(2^k-1)/2$, which changes all pairwise sums by
$-(2^k-1)$, the sets $E$ and $O$ are the sets of sums of an even,
respectively odd, number of the numbers $1, 2, 4, \ldots, 2^{k-1}$. So
$P(E) = P(O)$ by~\cite[Thm.~2]{selfstraus}, \cite[\S6]{gordonfs}, and since $E =
A(\lambda)$ and $O = A(\mu)$ we get $M(\lambda) = M(\mu)$. Since $E \cap O
= \emptyset$ we have $A(\lambda) \neq A(\mu)$, so $\mu \neq \lambda$ and
$\mu \neq \sigma(\lambda)$. For $l \geq 5$ the graph automorphism $\sigma$
is the only non-trivial one.
\end{Prf}

For $m = 1$, that is type $D_4$, we get $\lambda + \rho = 3\omega_1 +
2\omega_2 + \omega_3 + 4\omega_4$ and $\mu + \rho = \omega_1 + 2\omega_2 +
3\omega_3 + 4\omega_4$, which are related by a graph automorphism of order
$3$, as it must be by Theorem~\ref{thmClassical}. For $m = 2$, that is type
$D_{16}$, we get $\lambda + \rho$ and $\mu + \rho$ with coordinates
\begin{align*}
& (3,2,1,3,1,2,3,2,1,2,3,1,3,2,1,4) \quad \textrm{and} \\
& (1,2,3,1,3,2,1,2,3,2,1,3,1,2,3,4)
\end{align*}
with respect to the fundamental weights.

%

\section{AI disclosure}

The author has used the AI tool \emph{Claude Opus 5.5} by Anthropic
for preparing this manuscript.
Some major mathematical ideas where suggested by Claude: the idea to use
the known invariants for Weyl groups of type $B_l$, which lead to 
Lemma~\ref{laBC}; and the connections to combinatorics (homometric sets 
and Moser's problem) which lead to the examples in Remark~\ref{remAl} and
in Proposition~\ref{proDlCounter}.

Furthermore, Claude has drafted parts of the text from a (quite detailed) sketch, 
proofread  the final text, and always checked everything with independent computations.

Of course, the author takes full responsibility for the final text.

\end{document}